\documentclass[12pt,twoside,reqno]{amsart}

\usepackage{amsmath,amssymb,amsthm}
\usepackage{fancyhdr}
\usepackage{epsfig,color,colordvi} 

\usepackage[numbers,sort&compress,square,comma]{natbib}

\theoremstyle{definition}

\theoremstyle{remark}
\newtheorem{remark}{\sc Remark}

\theoremstyle{plain}
\newtheorem{theorem}{\sc Theorem}

\newtheorem{proposition}{\sc Proposition}
\newtheorem{corollary}{\sc Corollary}

\newcommand{\be}{\begin{equation}}
\newcommand{\ee}{\end{equation}}
\newcommand{\nn}{\nonumber}

\def\eps{\varepsilon}

\def\bbR{\mathbb R}

\def\bbZ{\mathbb Z}

\def\fQ{{{\rm Q}\kern-.65em {}^{{}_/ }\,}}
\def\fQQ{ {{\rm Q}\kern-.57em \scriptscriptstyle{}^{]\kern.055em[}\,}}

\def\ord{\kern0.1em o\kern-0.02em{}_{\ds\breve{}}\kern0.1em}
\def\Ord{\kern0.1em O\kern-0.02em{\ds\breve{}}\kern0.1em}
\def\ds{\displaystyle}

\def\fmonth{\ifcase\month\or Jan\or Feb\or Mar\or Apr
\or May\or Jun\or Jul\or Aug\or Sep
\or Oct\or Nov\or Dec\fi\ }
\def\mmddyyyy{\the\month.\the\day.\the\year}
\def\ddmmyyyy{\the\day.\the\month.\the\year}
\def\Mddyyyy{\fmonth~\the\day,~\the\year}

\def\R{\bbR}

\def\Z{\bbZ}

\providecommand{\abs}[1]{\left\vert#1\right\vert}
\providecommand{\norm}[1]{\left\Vert#1\right\Vert}

\providecommand{\Norm}[1]{\muskip0=-2mu{\left|\mkern\muskip0\left|
\mkern\muskip0\left|#1\right|\mkern\muskip0
\right|\mkern\muskip0\right|}}

\numberwithin{equation}{section}

\allowdisplaybreaks[1]

\begin{document}


\author{F.~Rassoul-Agha}   
\address{Mathematical Biosciences Institute, Ohio State University, 
Columbus, OH 43210}
\email{firas@math.ohio-state.edu}
\urladdr{www.math.ohio-state.edu/$\sim$firas}
\author{T.~Sepp\"al\"ainen}  
\address{Mathematics Department, University of Wisconsin-Madison, Madison, WI 53706}
\email{seppalai@math.wisc.edu}
\urladdr{www.math.wisc.edu/$\sim$seppalai}
\thanks{T.~Sepp\"al\"ainen was partially supported by
National Science Foundation grant DMS-0402231.}

\date{November 26, 2004}





\begin{abstract}
We prove an invariance principle for a vector-valued
additive functional of a Markov chain for almost every starting point
with respect to an ergodic equilibrium distribution.
The hypothesis is a moment bound on the resolvent.
\end{abstract}




\title[Invariance principle for Markov chains]
{An Almost Sure Invariance Principle for Additive Functionals of Markov Chains}

\maketitle


\def\X{\mathcal X}
\def\B{\mathcal B}

\section{Introduction.}
\label{intro}
This note extends a result of Maxwell and Woodroofe \cite{MW}.
Our notation and
presentation follow  \cite{MW} as closely as possible,
and some results  from there
will be repeated without proofs.
The work presented here was motivated
by applications to random walk in random environment that are reported 
elsewhere.

After completing this note we learned of the work of Derriennic and Lin 
on fractional coboundaries of Banach space contractions \cite{DL2}.
The estimates needed for the invariance principles we prove
can be then obtained by applying the Derriennic and Lin machinery, and this way 
one can even improve the moment hypothesis to just having two moments 
($p=2$ below); see \cite{DL}.  
Thus, currently our note offers alternative probabilistic proofs of the 
results of \cite{DL} under the more restrictive moment hypothesis of \cite{MW}.

Let $(X_n)_{n\geq0}$ be a stationary  ergodic Markov chain
defined on a probability space $(\Omega,\mathcal{F},P)$,
with values
in a general measurable space $(\X,\B)$. Let $Q(x;dy)$ be
its transition probability kernel and
$\pi$ the stationary marginal distribution of each $X_n$.
Write $E$ for the expectation under $P$.
 $P_x$ denotes the probability measure obtained by
conditioning on $X_0=x$, and
 $E_x$ is the corresponding expectation. For $p\geq1$, we will
denote by $L^p(\pi)$ the equivalence class
of $\B$-measurable functions $g$ with values
in $\R^d$ for some $d\geq1$ and such that
$$\norm{g}_p^p=\int|g(x)|^p\pi(dx)<\infty.$$
Here, $|\cdot|$ denotes the $\ell^2$-norm on $\R^d$.

Now fix $d$ and an $\R^d$-valued   function $g\in L^2(\pi)$ with
$\int g\,d\pi=0$. Define $S_0(g)=0$ and
$$S_{n+1}(g)=\sum_{k=0}^{n}g(X_k)
\hbox{ and }\widetilde S_n(g)=S_n(g)-E_{X_0}(S_n(g)),
~\hbox{for }n\geq0.$$

We are concerned with central limit type results
for $S_n(g)$ and $\widetilde S_n(g)$. 
This question has been investigated from many
angles and under different assumptions; see \cite{MW} and
its references.  A widely used method of
Kipnis and Varadhan \cite{KV}
works for reversible chains.
Article \cite{MW}  adapted this approach to a
non-reversible setting, and used
growth bounds on the resolvent to  obtain
 sufficient conditions  for an invariance principle for $S_n(g)$ under $P$,
if $p>2$.

Derriennic and Lin used then their theory of fractional coboundaries \cite{DL2}
to push the result to an invariance principle for $S_n(g)$ under $P_x$, for 
$\pi$-a.e.\ $x$, even when $p=2$; see \cite{DL}. Using their method one can 
also show that the same almost-sure invariance principle holds for
$\widetilde S_n(g)$.
We will show how to further the probabilistic technique of \cite{MW} 
to yield both almost-sure invariance principles 
(for $S_n(g)$ and $\widetilde S_n(g)$) when $p>2$.

Invariance principles
for additive functionals of Markov chains
 have many applications. This note is a byproduct
of the authors' recent work on random walks in a random environment
\cite{qclt-spacetime,forbidden} where this invariance principle proved useful.

Let us now describe the structure of this note. In Section \ref{martingales}
we will present the setting of \cite{MW} and prove an $L^q$ bound,
with $q>2$, on a certain martingale. In Section \ref{clt} we will state and
prove the main theorem of the note. The proof depends on
 a vector-valued version of a well-known
invariance principle for martingales (Theorem 3 of \cite{qclt-spacetime}).  

\section{A useful martingale.}
\label{martingales}
For a function $h\in L^1(\pi)$ and $\pi$-a.e.~$x\in\X$ define
$$Qh(x)=\int h(y)Q(x;dy).$$
$Q$ is a contraction on $L^p(\pi)$ for every $p\geq1$.
For $\eps>0$ let $h_\eps$ be the solution of
$$(1+\eps)h_\eps-Qh_\eps=g.$$
In other words,
$$h_\eps=\sum_{k=1}^\infty(1+\eps)^{-k}Q^{k-1}g.$$
Note that $h_\eps\in L^p(\pi)$, if $g\in L^p(\pi)$.
On $\X^2$ define the function
$$H_\eps(x_0,x_1)=h_\eps(x_1)-Qh_\eps(x_0).$$
For a given realization of $(X_k)_{k\geq0}$, let
$$M_n(\eps)=\sum_{k=0}^{n-1}H_\eps(X_k,X_{k+1})\ \hbox{ and }\
R_n(\eps)=Qh_\eps(X_0)-Qh_\eps(X_n)$$
so that
$$S_n(g)=M_n(\eps)+\eps S_n(h_\eps)+R_n(\eps).$$
Finally, let $\pi_1$ be the distribution of $(X_0,X_1)$ under $P$; that is
$$\pi_1(dx_0,dx_1)=Q(x_0;dx_1)\pi(dx_0).$$ 
Let us denote the $L^p$-norm on
$L^p(\pi_1)$ by $\Norm{\cdot}_p$.
The following theorem summarizes results of \cite{MW}.

\newtheorem*{theoremMW}{Theorem MW}
\begin{theoremMW}
Assume that $g\in L^2(\pi)$ and that there exists an $\alpha\in(0,1/2)$ such 
that
\begin{align}
\label{cond}
\norm{\sum_{k=0}^{n-1}Q^k g}_2=\Ord(n^\alpha).
\end{align}
Then
\begin{enumerate}
\item The limit $H=\lim_{\eps\rightarrow0^+}H_\eps$ exists in
$L^2(\pi_1)$. Moreover, if one defines
$$M_n=\sum_{k=0}^{n-1}H(X_k,X_{k+1}),$$
then, for $\pi$-almost every $x$,  $(M_n)_{n\geq1}$ is a $P_x$-square
integrable martingale, relative to the filtration
$\{{\mathcal F}_n=\sigma(X_0,\cdots,X_n)\}_{n\geq0}$.
\item One has $\norm{h_\eps}_2=\Ord(\eps^{-\alpha})$, and
if $R_n=S_n(g)-M_n=M_n(\eps)-M_n+\eps S_n(h_\eps)+R_n(\eps)$, then
$$E(\abs{R_n}^2)=\Ord(n^{2\alpha}).$$
\end{enumerate}
\end{theoremMW}

\begin{proof}
The existence of $H$ follows from Proposition 1 of \cite{MW}.
The statement about $M_n$ follows from Theorem 1 therein.
The bounds on $\norm{h_\eps}_2$ and $E(\abs{R_n}^2)$ follow from Lemma 1
and Corollary 4 of \cite{MW}, respectively.
%
%
\end{proof}

If, moreover, one has an $L^p$ assumption on $g$, then one can say more.

\begin{theorem}
\label{Lq}
Assume that there exists an $\alpha<1/2$ for which {\rm(\ref{cond})}
is satisfied.
Assume also that there exists a $p>2$ such that $g\in L^p(\pi)$. Then
there exists a $q\in(2,p)$ such that $H\in L^q(\pi_1)$ and
$(M_n)_{n\geq1}$ is an $L^q$-martingale.
\end{theorem}

\begin{proof}
First choose a positive $q<(3-2\alpha)p/(1-2\alpha+p)$.
One can check that since $2\alpha<1$ and $p>2$, we have $q\in(2,p)$.
Using H\"older's inequality, we have
$$\Norm{H_\delta-H_\eps}_q^q\leq\Norm{H_\delta-H_\eps}_p^a~
\Norm{H_\delta-H_\eps}_2^b,$$
where $a=p(q-2)/(p-2)<q$ and $b=q-a$.
Next, observe that
$$\norm{h_\eps}_p\leq\sum_{n\geq1}(1+\eps)^{-n}\norm{g}_p
=\norm{g}_p\eps^{-1}.$$
Thus, one has
$$\Norm{H_\delta-H_\eps}_q^q\leq2^a\norm{g}_p^a(\eps^{-1}+\delta^{-1})^a
\Norm{H_\delta-H_\eps}_2^b,$$
and, by Lemma 2 of \cite{MW},
$$\Norm{H_{\delta_k}-H_{\delta_{k-1}}}_q^q\leq C\,
2^{ka}\cdot 2^{-kb/2}
(\norm{h_{\delta_k}}^2_2+\norm{h_{\delta_{k-1}}}_2^2)^{b/2},$$
where $\delta_k=2^{-k}$. By part (ii) of Theorem MW, we know that
$\norm{h_\delta}_2=\Ord(\delta^{-\alpha})$. Therefore, one has
$$\Norm{H_{\delta_k}-H_{\delta_{k-1}}}_q^q\leq C\,2^{k(a-b/2+\alpha b)},$$
with maybe a different $C$ than above. Now, by the choice of $q$,
one can verify that $a-b/2+\alpha b<0$, and then
repeat the proof of Proposition 1 in \cite{MW}, with $\Norm{\cdot}_2$
replaced by $\Norm{\cdot}_q$.
\end{proof}

\begin{remark}
Note that {\rm \cite{MW}} uses $\norm{\cdot}_1$ for
the $L^2$-norm under $\pi_1$,
while we use $\Norm{\cdot}_2$.
\end{remark}

\section{The almost sure invariance principle.}
\label{clt}
First some notation. We write $A^T$ for the transpose of a vector or matrix
$A$. An element of $\R^d$ is regarded as a $d\times 1$ matrix,
or column vector.
Define
$${\mathbb B}_n(t)=n^{-1/2}S_{[nt]}(g)\hbox{ and }
\widetilde{\mathbb B}_n(t)=n^{-1/2}\widetilde S_{[nt]}(g),\hbox{ for }t\in[0,1].$$
Here, $[x]=\max\{k\in\Z:k\leq x\}$.
Let $D_{\R^d}([0,1])$ denote the space of
right continuous functions on
$[0,1]$ taking values in $\R^d$ and having left limits. 
This space is endowed with
the usual Skorohod topology \cite{billingsley}.
Let $\Delta$ denote  the Prohorov metric on
the space of Borel probability measures on $D_{\R^d}([0,1])$.

For a given symmetric, non-negative definite $d\times d$ matrix $\Gamma$,
a Brownian motion with diffusion matrix
$\Gamma$ is the $\R^d$-valued process $\{W(t):0\le t\le1\}$
such that $W(0)=0$, $W$ has continuous paths, independent increments,
and for $s<t$ the $d$-vector $W(t)-W(s)$ has Gaussian distribution
with mean zero and covariance matrix
$(t-s)\Gamma$. If the rank of $\Gamma$ is $m$,
one can produce such a process by
finding a $d\times m$ matrix $\Lambda$ such that
$\Gamma=\Lambda\Lambda^T$, and by defining $W(t)=\Lambda B(t)$
where $B$ is an  $m$-dimensional standard Brownian motion.

Let $\Phi_\Gamma$ denote
the distribution of Brownian motion with diffusion matrix $\Gamma$
on the space $D_{\R^d}([0,1])$. For $x\in\X$ let
$\Psi_n(x)$, respectively $\widetilde\Psi_n(x)$,
be the distribution of ${\mathbb B}_n$, respectively $\widetilde{\mathbb B}_n$,
on the Borel sets of $D_{\R^d}([0,1])$ under the measure $P_x$;
that is, conditioned on $X_0=x$.

Here is our main  theorem.

\begin{theorem}
\label{MW}
Assume
there are $p>2$ and $\alpha<1/2$ for which $g\in L^p(\pi)$ and
$E(|R_n^2|)=\Ord(n^{2\alpha})$.
Then
$$\lim_{n\rightarrow\infty}\Delta(\Phi_{\mathfrak D},\Psi_n(x))=0 \
\hbox{for }\pi\hbox{-a.e.\ }x,$$
where ${\mathfrak D}=E(M_1M_1^T)=\int HH^T\,d\pi_1.$
\end{theorem}

\begin{remark}
\label{improve}
The above result  improves  Theorem {\rm2}
 of {\rm \cite{MW}} which stated
 that
$$\lim_{n\rightarrow\infty}\int\Delta(\Phi_{\mathfrak D},\Psi_n(x))\pi(dx)=0.$$
\end{remark}

\begin{remark}
Due to Theorem {\rm MW}, {\rm(\ref{cond})} guarantees the bound on
$E(|R_n|^2)$ in Theorem \ref{MW}.
\end{remark}

\begin{proof} The proof is essentially done in \cite{MW}.
We explain below how to apply Borel-Cantelli's Lemma
to strengthen their result to an almost sure
statement.

Let $M_n^*(t)=n^{-1/2}M_{[nt]}$. We have
$$
\sup_{0\le t\le 1} \abs{{\mathbb B}_n(t)-M_n^*(t)} \le
n^{-1/2}\max_{k\le n}\abs{R_k}.
$$
 Therefore to conclude the proof we need to show two things:
\begin{align}
&\mbox{for $\pi$-almost every $x$,
under the probability measure $P_x$ the processes  }\nonumber\\
&\mbox{$M_n^*$ converge weakly to a Brownian motion
with diffusion matrix $\mathfrak{D}$,}
\label{M-lim}
\end{align}
and
\begin{align}
\label{R}
n^{-1/2}\max_{k\le n}\abs{R_k}
\mathop{\longrightarrow}_{n\rightarrow\infty}0 \
\hbox{ in }P_x\hbox{-probability, for }\pi\hbox{-a.e.\ }x.
\end{align}

Statement (\ref{M-lim}) follows from the martingale invariance
principle  stated  as Theorem 3 in \cite{qclt-spacetime}.
The limits needed as hypotheses for that theorem follow from
ergodicity and the square-integrability of $H$. We leave
this check to the reader.

To prove (\ref{R}), let  $n_j=j^r$ for a large enough integer $r$.
Fix $0<\gamma<1$, and let
$m_j=\lceil n_j^{1-\gamma}\rceil$, $\ell_j=\lceil n_j^{\gamma}\rceil$.
Here $\lceil x\rceil=\min\{n\in\Z:x\leq n\}.$
Since $R_n=S_n(g)-M_n$, one can write
\begin{align}
n_j^{-1/2}\max_{i\le n_j}\abs{R_i}&\leq
n_j^{-1/2}\max_{0\le k\le m_j}\abs{R_{k\ell_j}}\nn\\
&\quad
+n_j^{-1/2}\max_{0\le k<m_{j}}\,\max_{k\ell_j\le i\le (k+1)\ell_j}
\abs{M_i-M_{k\ell_j}}\nn\\
&\quad
+n_j^{-1/2}\max_{0\le k<m_{j}}\,\max_{k\ell_j\le i\le (k+1)\ell_j}
\abs{S_i(g)-S_{k\ell_j}(g)}.
\label{Rn}
\end{align}

Recalling that $E(\abs{R_n}^2)=\Ord(n^{2\alpha})$ with $\alpha<1/2$,
one can
apply Corollary 3 of \cite{MW} to get that for any $\delta>0$
$$P(\max_{0\le k\le m_j}\abs{R_{k\ell_j}}\geq\delta\sqrt{n_j})=
\Ord(\ell_j^{2\alpha} m_j^\beta/n_j)=
\Ord(j^{-r(1-2\gamma\alpha-(1-\gamma)\beta)}),$$
for any $\beta>1$. Choosing $\beta$ close enough to $1$
 and $r$ large enough,
the above becomes summable.
Borel-Cantelli's Lemma implies then that the first term on the
right-hand-side of (\ref{Rn}) converges to $0$, $P$-a.s.

The second martingale term on the right-hand side
of (\ref{Rn}) tends to $0$ in $P_x$-probability
for $\pi$-a.e.~$x$, by the
functional central limit theorem for $L^2$-martingales;
see Theorem 3 of \cite{qclt-spacetime}, for example. 
So it all boils down to showing that
the last term in (\ref{Rn}) goes to $0$ $P$-a.s.

\begin{remark}
\label{where}
Note that we have so far used the fact that $g\in L^2(\pi)$. It is only
to control the third term in {\rm(\ref{Rn})} that we need a higher moment.
\end{remark}

Define, for $\delta>0$,
$$B_j^{'}=\{\max_{0\le k<m_{j}}\,\max_{k\ell_j\le i\le (k+1)\ell_j}
\abs{S_i(g)-S_{k\ell_j}(g)}\geq\delta\sqrt{n_j}\}.$$
Since $g\in L^p(\pi)$, one can write:
\begin{align*}
P(B_j^{'})&\leq P(\,\max_{i\leq n_j}|g(X_i)|\geq
\delta \sqrt {n_j}/\ell_j)\\
&\leq n_j \pi(\,|g|\geq\delta\sqrt{n_j}/\ell_j)
=\Ord(j^{-r(p/2-1-\gamma p)}).
\end{align*}
By choosing $\gamma$  small enough and $r$ large enough, one can make sure that
$P(B_j^{'})$ is summable. By Borel-Cantelli's Lemma,
the third term in (\ref{Rn}) converges to $0$, $P$-a.s.

Finally, note that if $n_{j-1}\le n\le n_j$, then
\begin{align}
\label{trick}
\max_{k\leq n}\frac{\abs{R_k}}{\sqrt n}
\leq
\left(\frac{j}{j-1}\right)^{r/2}
\max_{k\le n_j}\frac{\abs{R_k}}{\sqrt{n_j}},
\end{align}
and so (\ref{R}) follows.
\end{proof}


\begin{remark}
In the above proof we only needed the martingale term in
{\rm (\ref{Rn})}
to converge in $P_x$-probability. The $L^q$-bounds of
Theorem {\rm\ref{Lq}} imply that it actually goes to $0$ $P$-a.s.,
making
{\rm(\ref{R})} also true $P$-a.s. All this is of course under the assumptions
$p>2$ and {\rm(\ref{cond})} with $\alpha<1/2$. In \cite{DL} it is shown
that the same almost-sure convergence happens even when $p=2$.
\end{remark}

We also have a similar result for $\widetilde S_n(g)$:

\begin{theorem}
\label{centered}
Assume there are $p>2$ and $\alpha<1/2$ for which $g\in L^p(\pi)$ and
condition {\rm(\ref{cond})} is satisfied.
Then $n^{-1/2}\max_{k\leq n}|E_x(S_k(g))|$ converges to $0$ as 
$n$ goes to infinity for $\pi$-almost every $x$. Consequently, for 
$\pi$-almost every $x$, 
$$\lim_{n\to\infty}
n^{-1/2}\max_{k\le n}|S_k(g)-\widetilde S_k(g)|=0~
P_x\hbox{-almost surely,}$$
and, therefore,
$$\lim_{n\rightarrow\infty}\Delta(\Phi_{\mathfrak D},\widetilde\Psi_n(x))=0 
\hbox{ for }\pi\hbox{-a.e.\ }x.$$
The diffusion matrix ${\mathfrak D}$ is as defined in Theorem {\rm\ref{MW}}.
\end{theorem}

Before we start the proof,
we need to reprove a maximal inequality of \cite{MW},
this time for a Markov transition operator rather than a shift.
For a probability transition kernel $Q$ and a function $g$ in its domain, define
$T_n(g,Q)=\sum_{k=0}^{n-1}Q^k g.$
We then have the following:

\begin{proposition}
\label{maximal}
Let $Q$ be a probability transition kernel with invariant measure $\pi$.
Let $g\in L^2(\pi)$ be such that
$$\int|T_n(g,Q)|^2 d\pi\le C(g,Q)n,$$
for some $C(g,Q)<\infty$ and all $n\ge1$. Then we have
$$\pi\left(\max_{j\le n}|T_j(g,Q)|>\lambda\right)
\le\frac{2^{6k}C(g,Q) n^{1+2^{-k}}}{\lambda^2},$$
for all $n\geq1$, $k\geq0$, and $\lambda>0$.
\end{proposition}

\begin{proof}
We will proceed by induction on $k$. For $k=0$ the lemma follows from Chebyshev's
and Jensen's inequalities, as well as the invariance of $\pi$ under $Q$.
Let us assume that the lemma has been proved for some $k\geq0$. We will prove
it for $k+1$. To this end, choose $n\ge1$ and $\lambda>0$. Let
$m=\lceil{\sqrt n}\,\rceil$. Then $[n/m]\leq m$, and
\begin{align*}
\pi\left(\max_{j\le n}|T_j(g,Q)|>\lambda\right)
&\le \pi\left(\max_{i\le n/m}|T_{im}(g,Q)|>\lambda/2\right)\\
&\quad+m\max_{i\le n/m}
\pi\left(\max_{j\leq m}|T_{j+im}(g,Q)-T_{im}(g,Q)|>\lambda/2\right)\\
&\le \pi\left(\max_{i\le m}|T_i(T_m(g,Q),Q^m)|>\lambda/2\right)\\
&\quad+
m\max_{i\le m}\pi\left(\max_{j\leq m}|T_j(Q^{im}g,Q)|>\lambda/2\right)\\
&\leq\frac{4\cdot 2^{6k} C(T_m(g,Q),Q^m) m^{1+2^{-k}}}{\lambda^2}\\
&\quad+\max_{i\le m}\frac{4\cdot 2^{6k} C(Q^{im}g,Q) m^{2+2^{-k}}}{\lambda^2}.
\end{align*}
But one has
$$\int|T_n(T_m(g,Q),Q^m)|^2 d\pi=\int|T_{mn}(g,Q)|^2 d\pi\leq C(g,Q)mn$$
and, therefore, $C(T_m(g,Q),Q^m)\leq C(g,Q)m$. Similarly,
\begin{align*}
\int|T_n(Q^{im}g,Q)|^2 d\pi&=\int|Q^{im}T_n(g,Q)|^2 d\pi
\leq\int Q^{im}|T_n(g,Q)|^2 d\pi\\
&=\int|T_n(g,Q)|^2 d\pi\leq C(g,Q)n.
\end{align*}
Thus, $C(Q^{im}g,Q)\leq C(g,Q)$. Above, we have used Jensen's inequality
to bring $Q$ outside the square and then the fact that $\pi$ is invariant
under $Q$.
Now, we have
$$\pi\left(\max_{j\le n}|T_n(g,Q)|>\lambda\right) \leq
\frac{8\cdot 2^{6k} C(g,Q)m^{2+2^{-k}}}{\lambda^2}.$$
Since $m\leq 2\sqrt n$, it follows that
$$\pi\left(\max_{j\le n}|T_n(g,Q)|>\lambda\right) \leq
\frac{2^{6k+6}C(g,Q)n^{1+2^{-k-1}}}{\lambda^2}$$
which is the claim of the lemma, for $k+1$.
\end{proof}

The following is then immediate:

\begin{corollary}
For any $\beta>1$ there is a constant $\Gamma$, depending only on $\beta$, for which
$$\pi\left(\max_{j\le n}|T_j(g,Q)|>\lambda\right)\leq
\frac{\Gamma C(g,Q) n^\beta}{\lambda^2}$$
for all $\lambda>0$ and $n\geq1$.
\end{corollary}

We can now prove the theorem.

\begin{proof}[Proof of Theorem {\rm\ref{centered}}]
Observe that $E_x(S_n(g))=T_n(g,Q)$. Now, 
recall that $n_j=j^r$, for an integer $r$ 
large enough. Also, for $0<\gamma<1$ we
have $m_j=\lceil{n_j^{1-\gamma}}\rceil$, $\ell_j=\lceil{n_j^\gamma}\rceil$.
Then,
\begin{align}
&n_j^{-1/2}\max_{i\le n_j}|E_x(S_i(g))|\le
n_j^{-1/2}\max_{k\le m_j}|E_x(S_{k\ell_j}(g))|
\label{line1}\\
&\qquad\qquad+
n_j^{-1/2}\max_{k<m_j}\max_{k\ell_j\leq i\leq(k+1)\ell_j}
|E_x(S_i(g))-E_x(S_{k\ell_j}(g))|.
\label{line2}
\end{align}
For the first term, we can use the above corollary to write
\begin{align*}
\pi\left(n_j^{-1/2}\max_{k\le m_j}|E_x(S_{k\ell_j}(g))|>\eps\right)
&=\pi\left(n_j^{-1/2}\max_{k\le m_j}|T_{k\ell_j}(g,Q)|>\eps\right)\\
&=\pi\left(\max_{k\le m_j}|T_k(T_{\ell_j}(g,Q),Q^{\ell_j})|>
\eps\sqrt{n_j}\right)\\
&\le\frac{\Gamma C(T_{\ell_j}(g,Q),Q^{\ell_j})m_j^\beta}{\eps^2 n_j}\\
&\le\frac{\Gamma C}{\eps^2}\frac{\ell_j^{2\alpha}m_j^\beta}{n_j}=
\Ord(j^{-r(1-2\alpha\gamma-(1-\gamma)\beta)})
\end{align*}
since
$$\int|T_n(T_{\ell_j}(g,Q),Q^{\ell_j})|^2 d\pi=
\int|T_{n\ell_j}(g,Q)|^2 d\pi\leq C({n\ell_j})^{2\alpha}\leq C{\ell_j}^{2\alpha}n,$$
by (\ref{cond}).
If one chooses $\beta$ small enough and $r$ large enough, then the
term on line (\ref{line1}) goes to $0$, $\pi$-a.s., by Borel-Cantelli's Lemma.
For the term on line (\ref{line2}) we have
\begin{align*}
&\pi\left(n_j^{-1/2}\max_{k<m_j}\max_{k\ell_j\leq i\leq(k+1)\ell_j}
|E_x(S_i(g))-E_x(S_{k\ell_j}(g))|\geq\eps\right)\\
&\qquad\leq
\pi\left(\max_{i\leq n_j}|Q^i g|\geq \eps\sqrt{n_j}/\ell_j\right)\\
&\qquad
\leq n_j \max_{i\leq n_j}\pi\left(|Q^i g|\geq\eps\sqrt{n_j}/\ell_j\right)\\
&\qquad\leq\Ord(j^{-r(p/2-1-\gamma p)}),
\end{align*}
since
$$\int |Q^i g|^p d\pi\leq\int Q^i(|g|^p)d\pi=
\int|g|^p d\pi<\infty.$$
Using Borel-Cantelli's Lemma, we get that the term on line (\ref{line2}) also
converges to $0$, $\pi$-a.s., if one chooses $\gamma$ small enough and $r$ large
enough.

Therefore, we have shown that
$n_j^{-1/2}\max_{i\le n_j}|E_x(S_i(g))|$
converges to $0$, $\pi$-a.s. The claim of the theorem follows then as in
(\ref{trick}), by considering $n_j\leq n\le n_{j+1}$.
\end{proof}

\bibliographystyle{aop} 
\bibliography{refs}

\begin{thebibliography}{7}
\expandafter\ifx\csname natexlab\endcsname\relax\def\natexlab#1{#1}\fi
\expandafter\ifx\csname url\endcsname\relax
  \def\url#1{\texttt{#1}}\fi
\expandafter\ifx\csname urlprefix\endcsname\relax\def\urlprefix{URL }\fi

\bibitem{billingsley}
\textsc{Billingsley, P.} (1999).
\newblock \textit{Convergence of probability measures}.
\newblock 2nd edn. John Wiley \& Sons Inc., New York.

\bibitem{DL2}
\textsc{Derriennic, Y.} and \textsc{Lin, M.} (2001).
\newblock Fractional {P}oisson equations and ergodic theorems for fractional
  coboundaries.
\newblock \textit{Israel J. Math.} \textbf{123} 93--130.

\bibitem{DL}
\textsc{Derriennic, Y.} and \textsc{Lin, M.} (2003).
\newblock The central limit theorem for {M}arkov chains started at a point.
\newblock \textit{Probab. Theory Related Fields} \textbf{125} 73--76.

\bibitem{KV}
\textsc{Kipnis, C.} and \textsc{Varadhan, S. R.~S.} (1986).
\newblock Central limit theorem for additive functionals of reversible {M}arkov
  processes and applications to simple exclusions.
\newblock \textit{Comm. Math. Phys.} \textbf{104} 1--19.

\bibitem{MW}
\textsc{Maxwell, M.} and \textsc{Woodroofe, M.} (2000).
\newblock Central limit theorems for additive functionals of {M}arkov chains.
\newblock \textit{Ann. Probab.} \textbf{28} 713--724.

\bibitem{qclt-spacetime}
\textsc{Rassoul-Agha, F.} and \textsc{Sepp\"al\"ainen, T.}
  (2004{\natexlab{a}}).
\newblock An almost sure invariance principle for random walks in a space-time
  i.i.d. random environment.
\newblock \textit{Probab. Theory Related Fields} To appear.

\bibitem{forbidden}
\textsc{Rassoul-Agha, F.} and \textsc{Sepp\"al\"ainen, T.}
  (2004{\natexlab{b}}).
\newblock Ballistic random walk in a random environment with a forbidden
  direction.
\newblock \textit{Ann. Probab.} Submitted.

\end{thebibliography}




%
%



\end{document}